# ON THE INVERSE FIRST-PASSAGE-TIME PROBLEM FOR A WIENER PROCESS[1,2]

By Cristina Zucca and Laura Sacerdote

*University of Torino*

The inverse first-passage problem for a Wiener process $(W_t)_{t \geq 0}$ seeks to determine a function $b : \mathbb{R}_+ \to \mathbb{R}$ such that

$$\tau = \inf\{t > 0 | W_t \geq b(t)\}$$

has a given law. In this paper two methods for approximating the unknown function $b$ are presented. The errors of the two methods are studied. A set of examples illustrates the methods. Possible applications are enlighted.

**1. Introduction.** Many phenomena can be modeled as first passage time of suitable Markov processes through constant or time varying boundaries. The first-passage problem has a long history and a large number of applications that range from finance to engineering and biology [see, e.g., Ricciardi and Sato (1990) for some references]. Yet explicit solutions to the first-passage problem are known only in a limited number of special cases (including linear or quadratic boundaries).

In modeling generally one describes the dynamics of the involved variables via a suitable stochastic process $\{X_t, t \geq 0\}$ constrained by an assigned boundary $b : \langle 0, \infty \rangle \to \mathbb{R}$ satisfying $b(0+) \geq 0$ and investigates distribution features of the first passage time (FPT)

$$(1.1) \qquad \tau_b = \inf\{t > 0 | W_t \geq b(t)\}$$

Received November 2006; revised January 2008.

[1]Supported in part by PRIN-Cofin 2005.

[2]This paper is published despite the effects of the Italian law 133/08 (http://groups.google.it/group/scienceaction). This law drastically reduces public funds to public Italian universities, which is particularly dangerous for scientific free research, and it will prevent young researchers from getting a position, either temporary or tenured, in Italy. The authors are protesting against this law to obtain its cancellation.

of $X_t$ over $b$. This is the direct FPT problem. However, there are also instances when the underlying stochastic process is assigned, one knows or estimates the FPT distribution $F_b$ and wishes to determine the corresponding boundary shape. This is the inverse first passage time (IFPT) problem. Different applications of the IFPT problem can be listed. In a reliability context one could compare performances of alternative objects, characterized by their lifetime distribution comparing the shapes of the boundaries. In this case one should assume the same underlying stochastic process to describe the evolution toward the crash of the object and to identify the crash time as an observed FPT of the process through an unknown boundary. Other interesting applications can arise when the resulting boundary exhibits a periodic behavior. This fact can suggest particular features of the modeled phenomenon in a finance or in a neurobiological context. Despite its importance in applications, the literature focuses only on specific problems [cf. Sacerdote and Zucca ([2003a](#), [2003b](#)) and Sacerdote, Villa and Zucca ([2006](#))] while there is a lack of mathematical results. Here we try to cover this gap focusing on the inverse first passage time problem for the Wiener process.

The IFPT problem was brought to our attention by Professor Goran Peskir, who presented us its first formulation by A. Shiryaev in 1976 (during a Banach center meeting). The original Shiryaev question considered the case when $F_b$ is an exponential distribution. There exist two early papers dealing with the existence problem written by Dudley and Gutmann ([1977](#)) and Anulova ([1980](#)). These papers provide the existence of *some* stopping times for given $F_b$; however, these stopping times are not of the form ([1.1](#)) for some *function* $b$. Furthermore, in Capocelli and Ricciardi ([1972](#)) the properties characterizing the FPT distribution in the case of the constant boundary are discussed.

The problem of the existence and uniqueness of $b$ is still open. Through the paper we do not enter into the question of the existence of $b$ but assume that such a boundary exists, is unique and sufficiently regular (continuous or $C^1$).

The main aim of the present paper is to propose two algorithms for approximating the unknown function $b$ when $F_b$ is given. After some preliminaries presented in Section [2](#), in Section [3](#) we introduce a first approach to the IFPT problem. It is based on the idea to approximate $b$ by a piecewise-linear boundary $c$ for which the law of $\tau_c$ can be computed. It leads to a Monte Carlo method for determining $c$. The second approach (Section [5](#)) is based on the classic idea due to Volterra (around 1896) on how to approximate an integral equation (of the first kind) by a system of $n$ equations in $n$ unknowns. Since the integral equation for $b$ is nonlinear, the resulting system is also nonlinear. We propose a numerical method for approximating $b$ at finitely many points.

Both methods produce an approximation of the boundary value. Hence it is necessary to evaluate the respective errors, that is, the difference between the exact unknown value and the computed approximation. We limit the study of the error at the discretization time knots obtained with constant discretization step $h$. The error of the first method is discussed in Section 4. It is due to multiple causes but it is dominated by the error due to the substitution of the continuous boundary $b$ with a piecewise-linear boundary $c$ and it is $O(h^2)$. Furthermore, due to the use of the Monte Carlo method to evaluate some involved integrals, the resulting boundary is estimated with a fixed confidence level $\alpha$. The error of the second approach is evaluated in Section 6 and it is $O(h)$ or $O(h^2)$ depending upon the numerical method used to evaluate the involved integral. In Section 7 we illustrate the goodness of the proposed methods by means of a set of examples. Finally, in Section 8 we consider the boundary corresponding to the exponential FPT distribution and we give a numerical answer to the 1976 Shiryaev question.

**2. First passage time: the direct and inverse problem.** Given a standard Wiener process $W = (W_t)_{t\geq 0}$ started at zero, and a sufficiently regular function $b\colon \langle 0,\infty\rangle \to \mathbb{R}$ satisfying $b(0+) \geq 0$, denote

$$(2.1) \qquad \tau_b = \inf\{t > 0 | W_t \geq b(t)\}$$

the FPT of $W$ over $b$, and set

$$(2.2) \qquad f_b(t) = \frac{d}{dt}\mathbb{P}(\tau_b \leq t)$$

to denote its density function for $t > 0$. In the sequel we also need to consider the Wiener process starting in $x_0$ at time $t_0$; in this case we write $f_b(t|x_0, t_0)$ in spite of $f_b(t)$. Throughout we assume that all regularity assumptions ensuring the existence of the objects introduced and properties imposed are fulfilled.

By $p_{a,b}(t, x|s, y)$, we denote the transition probability density function of $W$ at $t$ constrained by the *absorbing* boundary $b$ over $[s, t]$ given that $W_s = y$, that is,

$$(2.3) \qquad p_{a,b}(t, x|s, y) = \frac{\partial}{\partial x}\mathbb{P}(W_t \leq x, \tau_b > t | W_s = y)$$

for $x \leq b(t)$ with $t > s \geq 0$ and $y < b(s)$ given and fixed.

Under the hypothesis that the stochastic process (in this case the Wiener process) is assigned, the *direct first-passage time problem* seeks to determine $F_b$ when $b$ is given. The *inverse first-passage time problem* seeks to determine $b$ when $F_b$ is given. It is interesting to note that the IFPT problem looks for the boundary $b(t)$ that is a deterministic time dependent function. Furthermore, as proved by Strassen (1967), if the boundary $b \in C^1(\mathbb{R}^+)$, then the probability distribution of the FPT is absolutely continuous with continuous density.

In the literature there exist some equations that relate quantities (2.2) and (2.3), allowing, in a limited number of cases, to solve the direct FPT problem [cf. Ricciardi and Sato (1990)]. In Section 2.1 we list some existing closed form results for the Wiener process. These results will be used in Sections 3 and 7 to numerically validate the reliability of the two algorithms proposed. In absence of the analytical solution, the direct problem can be solved numerically, making use of one of the algorithms proposed in the literature [Buonocore, Nobile and Ricciardi (1987) and Zucca (2002)]. In Section 7 we use the algorithm introduced in Buonocore, Nobile and Ricciardi (1987) to estimate the FPT p.d.f. for a set of assigned boundaries. These numerical evaluations make it possible to enlarge the test set for the validation of the IFPT algorithms proposed in this paper. Finally, we briefly recall in Section 2.2 some results that will be useful later on [cf. Peskir (2002)], since one of the two algorithms proposed for the IFPT method requires the knowledge of $b(0)$.

In the sequel we will assume that $t_0 = x_0 = 0$ when this gives no loss of generality.

2.1. *Known boundaries.* The first-passage time density $f_b$ for a Wiener process is known explicitly in a few cases. The following three will be of interest in the sequel:

1. *Linear boundary.* If the boundary is given by

$$c(t) = \alpha + \beta t \tag{2.4}$$

with $\alpha > x_0$ and $\beta \in \mathbb{R}$, then [cf. Doob (1949), page 397, and Malmquist (1954), page 526]

$$
\begin{aligned}
f_c(t|t_0, x_0) &= \frac{\alpha - x_0}{\sqrt{2\pi(t - t_0)^3}} e^{-(\alpha + \beta(t - t_0) - x_0)^2/(2(t - t_0))} \\
&= \frac{\alpha - x_0}{(t - t_0)^{3/2}} \varphi\left(\frac{\alpha + \beta(t - t_0) - x_0}{\sqrt{(t - t_0)}}\right)
\end{aligned}
\tag{2.5}
$$

for $t > 0$, where $\varphi(x) = (1/\sqrt{2\pi})e^{-x^2/2}$ is the standard normal probability density function. Note that (2.5) is known as the Bachelier–Lévy formula.

2. *Daniels boundary.* If the boundary is given by

$$d(t) = \frac{\alpha}{2} - \frac{t}{\alpha} \log\left(\frac{\beta}{2} + \sqrt{\frac{\beta^2}{4} + \gamma e^{-\alpha^2/t}}\right), \tag{2.6}$$

where $\alpha > 0$, $\beta \geq 0$ and $\gamma > -\beta^2/4$, then [cf. Daniels (1969)]

$$
\begin{aligned}
(2.7) \qquad f_d(t) &= \frac{1}{\sqrt{2\pi t^3}}\left(e^{-d(t)^2/(2t)} - \frac{\beta}{2}e^{-(d(t)-\alpha)^2/(2t)}\right) \\
&= \frac{1}{t^{3/2}}\left(\varphi\left(\frac{d(t)}{\sqrt{t}}\right) - \frac{\beta}{2}\varphi\left(\frac{d(t)-\alpha}{\sqrt{t}}\right)\right)
\end{aligned}
$$

for $t > 0$.

3. *Piecewise-linear boundary.* If the boundary is given by

$$
(2.8) \qquad c(t) = \alpha_i + \beta_i t \qquad \text{for } t \in [t_{i-1}, t_i] \text{ and } i \geq 1,
$$

where $t_i = t_0 + ih$, $h > 0$ and $\alpha_i, \beta_i \in \mathbb{R}$ with $t_0 \geq 0$. Setting $\alpha_{i+1} = \alpha_i + \beta_i t_i$, we get that $t \mapsto c(t)$ is continuous on $[t_0, \infty\rangle$. Let us denote by $c_i = c(t_i)$ the knots of $c$ for $i \geq 0$.

The transition density function of $W$ in $x_1, x_2, \ldots, x_n$ at $t_1, t_2, \ldots, t_n$ constrained by the *absorbing* piecewise-linear boundary $c$ over $[t_0, t_n]$ given that $W_{t_0} = x_0 < \alpha_1$ is [recall (2.3) above]

$$
\begin{aligned}
(2.9) \qquad & p_{a,c}(t_1, x_1; t_2, x_2; \ldots; t_n, x_n | t_0, x_0) \\
&= \prod_{i=1}^{n} p_{a,c}(t_i, x_i | t_{i-1}, x_{i-1}) \\
&= \prod_{i=1}^{n}(1 - e^{-2(c_i - x_i)(c_{i-1} - x_{i-1})/(t_i - t_{i-1})}) \\
&\qquad \times \frac{1}{\sqrt{2\pi(t_i - t_{i-1})}}\exp\left(-\frac{(x_i - x_{i-1})^2}{2(t_i - t_{i-1})}\right)
\end{aligned}
$$

for $(x_1, x_2, \ldots, x_n) \in \mathbb{R}^n$ with $x_i \leq c_i$ for $1 \leq i \leq n$ and $x_0 < c_0$ where $t_0 < t_1 < t_2 < \cdots < t_n$ are given and fixed. This implies

$$
\begin{aligned}
(2.10) \qquad & \mathbb{P}(W_{t_1} \in C_1, \ldots, W_{t_n} \in C_n, \tau_c > t_n | W_{t_0} = x_0) \\
&= \int_{C_1} \cdots \int_{C_n} p_{a,c}(t_1, x_1; \ldots; t_n, x_n | t_0, x_0)\, dx_1 \cdots dx_n
\end{aligned}
$$

for any Borel set $C_i \subseteq \langle-\infty, c_i]$ with $1 \leq i \leq n$.

The identity (2.9) is a well-known fact in the case $n = 1$ [cf. Doob (1949), equation (4.2) and Section 5, or Durbin (1971), Lemma 1]. The case of general $n \geq 2$ [cf. Wang and Pötzelberger (1997)] follows readily by induction arguments using that $W$ has stationary independent increments.

Since there are only a few cases where the FPT density function is known in closed form, there has been a big effort in the past to find numerical approximation for it. One of the most used algorithms has been presented in Buonocore, Nobile and Ricciardi ([1987](#)). This method solves a Volterra integral equation of the second kind derived from the Fortet equation but characterized by a nonsingular kernel. In Section [7](#) we apply this method to determine numerically the FPT density function in the case of a time periodic oscillating boundary.

2.2. *Limits at zero.* One of the key issues in the numerical treatment of the inverse first-passage problem is to know $b(0+)$ in terms of $f_b(0+)$ and vice versa. In this section we display some known results on this relation. We will make use of these facts in Sections [7](#) and [8](#) below.

In the notation of Section [2](#), let $W = (W_t)_{t \geq 0}$ be a standard Wiener process started at zero, let $b: \langle 0, \infty \rangle \to \mathbb{R}$ be a continuous function satisfying $b(0+) \geq 0$, and let us assume that the first-passage time $\tau_b$ from ([2.1](#)) has a continuous density function $f_b$ having a limit $f_b(0+)$ in $[0, \infty]$. Then the following facts are known to hold [cf. Peskir ([2002](#))].

If there are $\varepsilon > 0$ and $\delta > 0$ such that

$$(2.11) \qquad b(t) \geq \sqrt{(2 + \varepsilon) t \log(1/t)}$$

for all $t \in \langle 0, \delta \rangle$, then $f_b(0+) = 0$. If there is $\delta > 0$ such that

$$(2.12) \qquad b(t) \leq \sqrt{2t \log(1/t)}$$

for all $t \in \langle 0, \delta \rangle$, then $f_b(0+) = +\infty$. Moreover, if we define a boundary $g$ by setting

$$(2.13) \qquad g(t) = \sqrt{2t \log(1/t) + t \log \log(1/t) + ct}$$

for $t \in \langle 0, \delta_c \rangle$ with $c \in \mathbb{R}$ given and fixed, then the limit $f_g(0+)$ exists and is given by

$$(2.14) \qquad f_g(0+) = \frac{e^{-c/2}}{\sqrt{4\pi}}.$$

Observe that $g$ from ([2.13](#)) locally at zero lies between the two functions appearing on the l.h.s. of ([2.11](#)) and ([2.12](#)) respectively (where $\varepsilon > 0$ may be as small as one likes).

**3. First approach: a piecewise linear approximation via Monte Carlo simulation.** In this section we face the IFPT problem by looking for a piecewise-linear approximation of the unknown boundary $b$ with an approach that can be considered analogous to the method used by Durbin ([1971](#)) for the direct FPT problem. Thus, we determine the piecewise-linear boundary $c$ ([2.8](#)) that approximates the exact boundary $b$ corresponding to the given first-passage density $f_b$.

First we consider the inverse problem when $f_b(0+) = 0$ and we assume that $b(0+) > 0$ is known. Then, using the notation of Section 2 (with $t_0 = 0$), we set $c_0 = b(0+)$ and, thus, $\alpha_1 = b(0+)$ as well. When $b(0+)$ is unknown we can guess its value, but in this case the error at the initial knot can dominate the final one. Since the resulting algorithm uses the Monte Carlo approach to estimate involved integrals, henceforth we call it the PLMC (Piecewise Linear Monte Carlo) method. In the following we consider a time discretization $t_i = t_0 + ih$, $i = 1, 2, \ldots$, where $h$ is a positive constant. The choice $h$ constant is made to simplify our notation, but the method can be easily extended to a nonconstant $h$.

The driving idea of our algorithmic approach is to equate the probability of the FPT of $W$ through $c$ in $\langle t_{i-1}, t_i \rangle$ and the analogous probability determined by the given first-passage density $f_b$ in the same interval. The resulting equations allow to determine the coefficients $\alpha_i$ and $\beta_i$ in (2.8) successively on the intervals $\langle t_{i-1}, t_i \rangle$ for $i \geq 1$. The algorithm can be divided in two successive steps.

STEP 1.  We determine the value of $\beta_1$ in (2.8) such that the probability of the first-passage of $W$ through $c$ in $(0, t_1]$ equals the probability of the first-passage of $W$ through $b$ in $(0, t_1]$, that is, we look for the value of $\beta$ such that

$$(3.1) \qquad \int_0^{t_1} f_c(t)\,dt = \int_0^{t_1} f_b(t)\,dt,$$

where $f_c(t)$ is given by (2.5).

STEP 2.  Given $\alpha_1, \ldots, \alpha_n$ and $\beta_1, \ldots, \beta_n$ with $n \geq 1$, we set $\alpha_{n+1} = \alpha_n + \beta_n t_n$ and we determine the value of $\beta_{n+1}$ such that the probability of the first-passage time of $W$ through $c$ in $(t_n, t_{n+1}]$ equals the probability of the first-passage time of $W$ through $b$ in $(t_n, t_{n+1}]$, that is, we look for the value of $\beta$ such that

$$
\begin{aligned}
(3.2) \qquad & \int_{t_n}^{t_{n+1}} \int_{-\infty}^{c_n} \cdots \int_{-\infty}^{c_1} \left( f_c(t|t_n, x_n) \prod_{i=1}^{n} p_{a,c}(t_i, x_i | t_{i-1}, x_{i-1}) \right) dx_1 \cdots dx_n \, dt \\
& = \int_{t_n}^{t_{n+1}} f_b(t)\,dt,
\end{aligned}
$$

where $f_c(t|t_n, x_n)$ is the density function of the first-passage time of $W$ over $\alpha_{n+1} + \beta_{n+1} t$ for $t > t_n$ given that $W_{t_n} = x_n$, that is, it is given by (2.5) with $\alpha := c_n$ and $\beta = \beta_{n+1}$, while $p_{a,c}(t_i, x_i; t_{i-1}, x_{i-1})$ is given by (2.9) for $1 \leq i \leq n$ (where $t_0 = x_0 = 0$).

REMARK 3.1. The product appearing on the l.h.s. of (3.2) involves the known functions (2.9) depending upon $\alpha_1, \ldots, \alpha_n$ and $\beta_1, \ldots, \beta_n$ determined in the preceding step. The unknown $\alpha_{n+1}$ and $\beta_{n+1}$ appear in (3.2) only within the function $f_c(t|t_n, x_n)$. Furthermore, $\alpha_{n+1}$ can (by continuity of $c$ at $t_n$) be expressed in terms of $\beta_{n+1}$ and the known $\alpha_n$ and $\beta_n$ as follows:

$$(3.3) \qquad \alpha_{n+1} = \alpha_n + (\beta_n - \beta_{n+1})t_n$$

so that we only need to determine $\beta_{n+1}$ from the equation (3.2).

Let us now detail the two steps of the algorithm. In order to solve the equation (3.1) and (3.2) for the unknown $\beta_1$ and $\beta_n$ respectively, it is necessary to compute the integrals involved.

Discretizing the integral of the l.h.s. in (3.1) by a rectangular method [cf. Atkinson (1989)], we obtain a nonlinear function of the unknown $\beta_1$ while the r.h.s. of (3.1) is easily computable by standard means. Various approximate methods can then be used to solve the resulting equation. Here we use the middle point method.

In the successive steps, mainly when $n$ becomes large, the multiple integrals of the l.h.s. of (3.2) cannot be handled by standard numerical methods. To handle the problem, we adapt the Monte Carlo method proposed by Wang and Pötzelberger (1997) to our case.

For this, note that, using expression (2.5) for $f_c(t|t_n, x_n)$ and (3.3), we get

$$
\begin{aligned}
H(\beta_{n+1}; x_n) &:= \int_{t_n}^{t_{n+1}} f_c(t|t_n, x_n)\, dt \\
&= \int_{t_n}^{t_{n+1}} \frac{c_n - x_n}{(t - t_n)^{3/2}} \varphi\left(\frac{c_n - x_n + \beta_{n+1}(t - t_n)}{\sqrt{t - t_n}}\right) dt \\
&= 1 - \Phi\left(\frac{\beta_{n+1}(t_{n+1} - t_n) + (c_n - x_n)}{\sqrt{t_{n+1} - t_n}}\right) \\
&\quad + e^{-2\beta_{n+1}(c_n - x_n)} \Phi\left(\frac{\beta_{n+1}(t_{n+1} - t_n) - (c_n - x_n)}{\sqrt{t_{n+1} - t_n}}\right) \\
&= 1 - \Phi\left(\frac{\beta_{n+1}(t_{n+1} - t_n) + \alpha_n + \beta_n t_n - x_n}{\sqrt{t_{n+1} - t_n}}\right) \\
&\quad + e^{-2\beta_{n+1}(\alpha_n + \beta_n t_n - x_n)} \\
&\quad \times \Phi\left(\frac{\beta_{n+1}(t_{n+1} - t_n) - \alpha_n - \beta_n t_n + x_n}{\sqrt{t_{n+1} - t_n}}\right)
\end{aligned}
$$

(3.4)

for each $n \geq 1$. Here $\Phi(x) = \int_{-\infty}^x \varphi(z)\, dz$ is the integral of the standard normal density. It follows by (2.9) and (3.4) that the equation (3.2) can be

rewritten as follows:

$$(3.5) \quad E\left( H(\beta_{n+1}; X_n) \prod_{i=1}^{n} I(X_i \leq c_i) e^{-2(c_i - X_i)(c_{i-1} - X_{i-1})/(t_i - t_{i-1})} \right)$$
$$= \int_{t_n}^{t_{n+1}} f_b(t) \, dt,$$

where $X_i \sim N(0, t_i - t_{i-1})$, $i = 1, \ldots, n$, are independent random variables (and $X_0 \equiv 0$). A Monte Carlo method can now be used to estimate the l.h.s. of (3.5).

In order to find an approximation of $\beta_{n+1}$, we can use, as in the Step 1 above, the iterative middle point procedure. Using then (3.3), one obtains $\alpha_{n+1}$ and thus determines $c$ on $(t_n, t_{n+1}]$.

REMARK 3.2. For each $n \geq 1$ given and fixed, the equation (3.5) defines a nonlinear function of the unknown parameter $\beta_{n+1}$ depending on it only through $H$ from (3.4). Since the l.h.s. of (3.5) is monotone in $\beta_{n+1}$, equation (3.5) admits a unique solution $\beta_{n+1} \in \mathbb{R}$.

REMARK 3.3. Finally, we modify the algorithm to include the case $f_b(0+) > 0$. In this case the equation (3.1) cannot be used on $(0, t_1]$ since the piecewise-linear boundary $c$ must satisfy $c(0+) > 0$, while the condition $f_b(0+) > 0$ implies $c(0+) = 0$ (cf. Section 2.2). However, we can use (2.14) to estimate the boundary in a neighborhood of zero. The use of this expression on the first discretization interval allows then the iteration of the previous algorithm on the successive intervals if we disregard the possible crossing happened in the neighborhood of zero.

**4. PLMC method error estimation.** The errors involved in the approximation of the boundary values via the PLMC method can be classified in three types according to the causes generating them. The first one is due to the piecewise linear approximation of the boundary, the second one is due to the Monte Carlo approximation of the integrals on the l.h.s. of (3.2) and the third one is due to the root method used to compute the zeros of (3.2). The last error can be disregarded since it depends on a tolerance that can be fixed in advance to make it negligible with respect to the other errors. In the following subsections we successively focus the study on the other two errors involved in the method. We explicitly underline how these two errors have different mathematical natures. Indeed, the first one, as the root error, is purely numerical, being determined by our approximation of the unknown function via a piecewise linear one. On the contrary, the second error has a statistical nature since it is related with the evaluation of an integral via a Monte Carlo method.

4.1. *Error due to the piecewise linear approximation.* The PLMC method estimates the boundary values through a piecewise linear approximation, but we limit the study of the error at the knots $t_n$, $n = 1, 2, \ldots$. Furthermore, at this step we assume the absence of other causes of error and we define *error of the method at the $n$th knot* the distance between the boundary $b$ and its approximated value $c$ on the $n$th knot $|\epsilon_n^{\mathrm{PLMC}}| = |b(t_n) - c(t_n)|$, $n = 1, 2, \ldots$. In order to gain an estimate of such error, we first consider in Lemma 4.1 the error due to the boundary linearization on a single step. In doing this analysis we hypothesize to know the true value of the boundary on the previous intervals. A second step considered in Lemma 4.2 studies the propagation of the error, that is, we admit an error $\delta$ in the estimation of the boundary value at node $n - 1$ and we evaluate its consequences on the next node $n$. Finally, in Theorem 4.3 we determine the global error of the PLMC Method.

LEMMA 4.1. *Let $(W_t)_{t \geq 0}$ be a Wiener process bounded by a monotone concave (or convex) boundary $b \in \mathcal{C}^2([0, \infty))$. If we approximate in $(0, t_n]$ the boundary $b(t)$ with the boundary $\widehat{c}_n(t)$, for $n = 1, 2, \ldots$, defined by*

$$(4.1) \qquad \widehat{c}_n(t) = \begin{cases} b(t), & t \in (0, t_{n-1}], \\ b(t_{n-1}) + \hat{\beta}_n(t - t_{n-1}), & t \in [t_{n-1}, t_n], \end{cases}$$

*the resulting error at the knot $t_n$, as $h \to 0$, is*

$$|b(t_n) - \widehat{c}(t_n)| \sim O(h^2).$$

PROOF. We limit ourself to the case of a concave boundary since the convex case can be dealt in a similar way. Moreover, we split the proof in two parts, first taking into account the first discretization step and later a generic step.

STEP 1. Let us consider the following three straight lines (cf. Figure 1) on $(0, h]$:

1. $B_1(t)$: the tangent to $b(t)$ in $t = 0$,

$$y = \alpha_1 + b'(0)t.$$

2. $c(t)$: the linear boundary determined via the PLMC method in the first discretization interval $(0, h]$,

$$y = \alpha_1 + \beta_1 t.$$

3. $B_2(t)$: the secant through $(0, \alpha_1)$ and $(h, b(h))$,

$$y = \alpha_1 + \frac{b(h) - \alpha_1}{h} t.$$

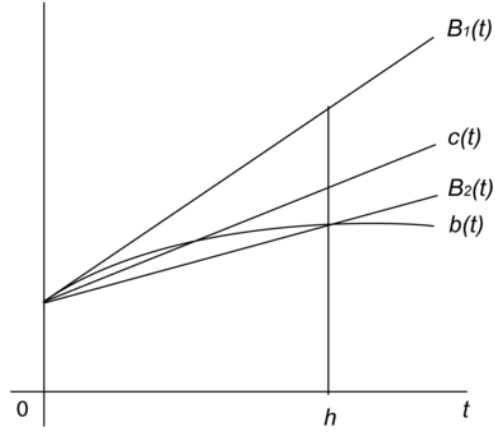

Fig. 1.  *Boundary $b(t)$, its tangent $B_1(t)$ in $t = 0$, the linear boundary $c(t)$ determined by the PLMC method and the secant $B_2(t)$ in $t = h$.*

Due to the concavity and monotonicity hypotheses on the boundary $b(t)$, it holds for all $t \in (0, h]$

$$B_2(t) \leq b(t) \leq B_1(t),$$

that implies [cf. Sacerdote and Smith (2004)]

$$\mathbb{P}(\tau_{B_1} \in (0, h]) \leq \mathbb{P}(\tau_c \in (0, h]) = \mathbb{P}(\tau_b \in (0, h]) \leq \mathbb{P}(\tau_{B_2} \in (0, h]).$$

Hence, we get

$$B_2(t) \leq c(t) \leq B_1(t),$$

that implies the ordering of the slopes

$$(4.2) \qquad \frac{b(h) - \alpha_1}{h} \leq \beta_1 \leq b'(0).$$

Note that $\beta_1$ is a function of $h$: $\beta_1(h)$. When $h \to 0$, due to the hypothesis $\alpha_1 = b(0)$, the inequality (4.2) implies

$$(4.3) \qquad \beta_1 = \beta_1(h) \xrightarrow[h \to 0]{} b'(0).$$

On the first knot the distance between $b(t)$ and the linear boundary of the PLMC method results in

$$|(\alpha_1 + \beta_1(h)h) - b(h)| = \left| (b(0) + \beta_1(h)h) - \left( b(0) + b'(0)h + \frac{b''(\xi)}{2!}h^2 \right) \right|$$

$$= \left| (\beta_1(h) - b'(0))h - \frac{b''(\xi)}{2!}h^2 \right|$$

$$\leq |(\beta_1(h) - b'(0))h| + \left| \frac{b''(\xi)}{2!}h^2 \right|$$

$$\leq |\beta_1'(\eta)h^2| + \left|\frac{b''(\xi)}{2!}h^2\right|$$

$$= O(h^2) + O(h^2) = O(h^2),$$

where we made use of (4.3), of the mean value theorem, of the triangular inequality and of the McLaurin expansion of $b(t)$ and $\beta_1(t)$. Here $\xi \in (0, h]$ and $\eta \in (0, h]$.

STEP $n$.    The boundary (4.1) coincides with the boundary $b(t)$ on $(0, t_{n-1})$ and it is determined solving

$$\mathbb{P}(\tau_b \in [t_{n-1}, t_n]) = \mathbb{P}(\tau_{\widehat{c}} \in [t_{n-1}, t_n]) \tag{4.4}$$

on $[t_{n-1}, t_n]$. Equation (4.4) implies

$$\begin{aligned}
&\mathbb{P}\{\tau_b \in [t_{n-1}, t_n] | X_t < b(t), t \leq t_{n-1}\} \mathbb{P}\{X_t < b(t), t \leq t_{n-1}\} \\
&\qquad = \mathbb{P}\{\tau_{\widehat{c}} \in [t_{n-1}, t_n] | X_t < \widehat{c}_n(t), t \leq t_{n-1}\} \mathbb{P}\{X_t < \widehat{c}_n(t), t \leq t_{n-1}\}
\end{aligned} \tag{4.5}$$

and, using (4.1), one has

$$\begin{aligned}
&\mathbb{P}(\tau_b \in [t_{n-1}, t_n] | X_t < b(t), t \leq t_{n-1}) \\
&\qquad = \mathbb{P}(\tau_{b(t_{n-1}) + \hat{\beta}_n(t - t_{n-1})} \in [t_{n-1}, t_n] | X_t < b(t), t \leq t_{n-1}).
\end{aligned} \tag{4.6}$$

The proof of Step 1 can now be repeated to the conditional probabilities in (4.6) in order to complete the proof.    □

Let us denote $\widehat{\alpha}_n$ the estimate of $\alpha_n$ and let $\widehat{\beta}_n(\widehat{\alpha}_n)$ be the corresponding estimate of the slope $\beta_n(\alpha_n)$ obtained with the PLMC method on the interval $[t_{n-1}, t_n]$.

LEMMA 4.2.    *An error* $|\delta| = |\widehat{\alpha}_n - \alpha_n|$ *propagates on the estimate* $\widehat{\beta}_n(\widehat{\alpha}_n)$:

$$|\Delta \beta_n| = |\beta_n(\alpha_n) - \widehat{\beta}_n(\widehat{\alpha}_n)| \sim O(\delta).$$

PROOF.    We separate the proof in two parts; the first one concerning the case of the first interval and the second one concerning the generic $n$th interval.

STEP 1.    Let us first prove that $\beta_1 = \beta_1(\alpha_1)$ is a continuous and differentiable function. When $\alpha_1 = b(0)$ define

$$F_1(\alpha_1, \beta_1) := 1 - \int_{-\infty}^{\alpha_1 + \beta_1 t_1} (1 - e^{-2\alpha_1/t_1(\alpha_1 + \beta_1 t_1 - x)}) \frac{1}{\sqrt{2\pi t_1}} e^{-x^2/(2t_1)} \, dx - k_1,$$

where $k_1$ is the r.h.s. of (3.1). Recognizing the first two terms as the l.h.s. of (3.1), equation (3.1) becomes $F_1(\alpha_1, \beta_1) = 0$. This is an implicit equation in $\beta_1$ and $\alpha_1$ that admits a continuous and differentiable solution $\beta_1 = \beta_1(t)$ if Dini's theorem holds. The hypothesis of Dini's theorem is verified since

$$\frac{\partial F_1(\alpha_1, \beta_1)}{\partial \beta_1} = -\int_{-\infty}^{\alpha_1 + \beta_1 t_1} \frac{2\alpha_1}{\sqrt{2\pi t_1}} e^{-2\alpha_1/t_1(\alpha_1 + \beta_1 t_1 - x) - x^2/(2t_1)} \, dx$$

is nonzero for all $\alpha_1, \beta_1$, having for the hypothesis $\alpha_1 \neq 0$. Hence, we get that $\beta_1 = \beta_1(\alpha_1)$ is a continuous and differentiable function and the error connected with the use of $\widehat{\alpha}_1$ in spite of $\alpha_1$ is propagated on $\beta_1$. Hence, $\beta_1$ becomes $\widehat{\beta}_1 = \beta_1 + \Delta\beta_1$ with

$$\Delta\beta_1 = \beta_1(\alpha_1 + \delta) - \beta_1(\alpha_1) \sim O(\delta)$$

since $\beta_1(\alpha_1)$ is continuous and differentiable.

STEP $n$. For $n = 2, 3, \ldots$, we proceed in analogy with step one showing that $\beta_n = \beta_n(\alpha_n)$ is a continuous and differentiable function. Let us consider the approximated stepwise linear boundary $c(t)$ in the time interval $[0, t_n]$. Let us define

$$F_n(\alpha_n, \beta_n)$$
$$:= \int_0^{c_1} \cdots \int_0^{c_{n-1}} dx_1 \cdots dx_{n-1} \prod_{i=0}^{n-1} f_a(x_i, t_i | x_{i-1}, t_{i-1})$$
$$\times \Big[1 - \int_{-\infty}^{\alpha_n + \beta_n t_n} (1 - \exp\{(-2(c_{n-1} - x_{n-1})$$
$$\times (\alpha_n + \beta_n t_n - x_n))$$
$$/(t_n - t_{n-1})\})$$
$$\times \exp\Big\{-\frac{(x_n - x_{n-1})^2}{(2(t_n - t_{n-1}))}\Big\} \Big/ (\sqrt{2\pi(t_n - t_{n-1})}) \, dx_n\Big]$$
$$- k_n,$$

so that equation (3.2) reads $F_n(\alpha_n, \beta_n) = 0$. This is an implicit equation in $\beta_n$ and $\alpha_n$ that admits a continuous and differentiable solution $\beta_n = \beta_n(t)$ if Dini's theorem holds.

Note that

$$\frac{\partial F_n(\alpha_n, \beta_n)}{\partial \beta_n}$$
$$= -\int_0^{c_1} \cdots \int_0^{c_{n-1}} \prod_{i=0}^{n-1} f_a(x_i, t_i | x_{i-1}, t_{i-1})$$
$$\times \int_{-\infty}^{\alpha_n + \beta_n t_n} 2t_n(c_{n-1} - x_{n-1})$$

$$\times\, e^{-2(c_{n-1}-x_{n-1})(\alpha_n+\beta_n t_n-x_n)/(t_n-t_{n-1})}$$

$$\times\, e^{-(x_n-x_{n-1})^2/(2(t_n-t_{n-1}))}$$

$$\Big/\big((t_n-t_{n-1})\sqrt{2\pi(t_n-t_{n-1})}\big)\, dx_n$$

is nonzero for all $\alpha_n, \beta_n$, since $x_{n-1} < c_{n-1}$. Hence, Dini's theorem holds and we get that $\beta_n = \beta_n(\alpha_n)$ is a continuous and differentiable function. If $\widehat{\alpha}_n = \alpha_n + \delta$, the error on $\alpha_n$ is propagated on $\beta_n$ as

$$\widehat{\beta}_n = \beta_n + \Delta\beta_n,$$

where

$$\Delta\beta_n = \beta_n(\alpha_n+\delta) - \beta_n(\alpha_n) \sim O(\delta). \qquad \square$$

THEOREM 4.3. *The error of the PLMC method at the discretization knots* $t_n$, $n = 1, 2, \ldots,$ *is*

$$|\epsilon_n^{\mathrm{PLMC}}| = |b(t_n) - c(t_n)| \sim O(\max(\delta, h^2))$$

*when* $\alpha_1$ *is affected by an error of the order of* $\delta \geq 0$.

PROOF. When $\alpha_1$ is affected by an error of the order of $\delta$, by Lemmas 4.1 and 4.2, we get that

$$|b(t_1) - c(t_1)| \leq |b(t_1) - \widehat{c}_1(t_1)| + |\widehat{c}_1(t_1) - c(t_1)|$$

$$\sim O(h^2) + O(\delta)$$

$$\sim O(\max(\delta, h^2)).$$

By induction, let $|b(t_{n-1}) - c(t_{n-1})| \sim O(\max(\delta, h^2))$, then

$$|\widehat{c}_n(t_n) - c(t_n)| = |b(t_{n-1}) + \widehat{\beta}_n(t_n - t_{n-1}) - c_{n-1} - \beta_n(t_n - t_{n-1})|$$

$$= |b(t_{n-1}) - c_{n-1} + (\widehat{\beta}_n - \beta_n)h|$$

$$\leq |b(t_{n-1}) - c_{n-1}| + |(\widehat{\beta}_n - \beta_n)h|$$

$$\sim O(\max(\delta, h^2)) + O(\max(\delta h, h^3))$$

$$\sim O(\max(\delta, h^2))$$

and we obtain

$$|b(t_n) - c(t_n)| \leq |b(t_n) - \widehat{c}_n(t_n)| + |\widehat{c}_n(t_n) - c(t_n)|$$

$$\sim O(h^2) + O(\max(\delta, h^2))$$

$$\sim O(\max(\delta, h^2)). \qquad \square$$

REMARK 4.1. If the boundary $b(t)$ is not monotone nor convex (or concave) but is sufficiently well behaved, one can proceed with a similar reasoning on each monotonicity and convexity interval. Hence, choosing $h$ in a suitable way, the stepwise boundary has still an error $|\epsilon_n^{\mathrm{PLMC}}| \sim O(\max(\delta, h^2))$.

4.2. *Error due to the Monte Carlo approximation.* Disregarding the error related to the numerical integration of the l.h.s. of (3.1) at the first step of the algorithm that can be well controlled with a careful use of numerical approximations, we focus on the next steps when the multiple integrals of the l.h.s. of (3.2),

$$
F(\beta_{n+1}) = \int_{t_n}^{t_{n+1}} \int_{-\infty}^{c_n} \cdots \int_{-\infty}^{c_1} \left( f_c(t|t_n, x_n) \right.
$$
$$
\left. \times \prod_{i=1}^{n} p_{a,c}(t_i, x_i|t_{i-1}, x_{i-1}) \right) dx_1 \cdots dx_n \, dt,
$$

are evaluated via Monte Carlo method. Using the Law of Large Numbers, we approximate the expectation on the l.h.s. of (3.5) with its sample mean. Hence, for a fixed confidence level $\alpha$, we get

$$
P\left( \left| \frac{1}{M} \sum_{j=1}^{M} H(\beta_{n+1}; X_{n,j}) \right. \right.
$$
$$
\times \prod_{i=1}^{n} I(X_{i,j} \leq c_i) e^{-2(c_i - X_{i,j})(c_{i-1} - X_{i-1,j})/(t_i - t_{i-1})}
$$
$$
\left. \left. - F(\beta_{n+1}) \right| < \delta_\alpha \right) > \alpha,
$$

where $X_{n,j} \sim N(0, t_i - t_{i-1})$, $i = 1, \ldots, n$, $j = 1, \ldots, M$, are independent random variables. Letting $x_{n,j}$, $i = 1, \ldots, n$, $j = 1, \ldots, M$, be a sample of $X_{n,j}$, with an accuracy $\delta_\alpha$ in the computation of the integral, we obtain a confidence interval for $\beta_{n+1}$ at the same level $\alpha$:

$$
\left[ F^{-1}\left( \frac{1}{M} \sum_{j=1}^{M} H(\beta_{n+1}; x_{n,j}) \right. \right.
$$
$$
\left. \left. \times \prod_{i=1}^{n} I(x_{i,j} \leq c_i) e^{-2(c_i - x_{i,j})(c_{i-1} - x_{i-1,j})/(t_i - t_{i-1})} \pm \delta_\alpha \right) \right].
$$

REMARK 4.2. The computations necessary to get the Monte Carlo evaluations are not excessively expensive. Hence, we can choose a very large value for the size $M$ for the number of simulations involved in the integral estimation. This allows to make this error negligible with respect to the error determined by the piecewise linear approximation of the boundary, but as a consequence of the use of the Monte Carlo method, our results are characterized by a reliability $\alpha$.

**5. An approximation by the nonlinear Volterra integral equation.** The algorithm discussed in the previous section is reliable and easily implemented, but it is computationally expensive since the Monte Carlo method requires long times of computation. In this section we therefore consider an alternative approach of pure numerical nature which is computationally more attractive. Since this method is based on the approximation of a non linear Volterra integral equation, it will be referred to as VIE (Volterra Integral Equation) method.

Let us consider the integral equation [cf. Peskir (2002)]:

$$(5.1) \qquad \Psi\left(\frac{b(t)}{\sqrt{t}}\right) = \int_0^t \Psi\left(\frac{b(t) - b(s)}{\sqrt{t-s}}\right) f_b(s)\,ds \qquad (t > 0),$$

where $\Psi(x) = 1 - \Phi(x)$ and $\Phi(x) = \int_{-\infty}^x \varphi(z)\,dz$ is the standard normal distribution. Equation (5.1) is a Volterra integral equation of the first kind in $f_b$, but it is a *nonlinear* Volterra integral equation of the second kind in $b$ and its kernel is nonsingular in the sense that it is bounded. Moreover the nonlinear functions involved in the equation are bounded and $\Psi$ is invertible. These features allow the development of the following numerical algorithm that approximate (5.1) in a simple and reliable way.

We numerically solve this equation to evaluate the approximate value $b^*$ of $b$ at the knots $t_i = ih$ for $i = 1, \ldots, n$, where $h = t/n$ (and $t > 0$ is given and fixed). To this aim we follow the original idea of Volterra [see, e.g., Linz (1985), Chapter VII], that is, we approximate the integral on the l.h.s. of (5.1) with the Euler method

$$(5.2) \quad \Psi\left(\frac{b^*(t_i)}{\sqrt{t_i}}\right) = \sum_{j=1}^i \Psi\left(\frac{b^*(t_i) - b^*(t_j)}{\sqrt{t_i - t_j}}\right) f_b(t_j) h \qquad (i = 1, \ldots, n),$$

getting a nonlinear system of $n$ equations in $n$ unknowns $b^*(t_1), \ldots, b^*(t_n)$. The $i$th equation of (5.2) for $i \geq 2$ makes use of the values $b^*(t_1), \ldots, b^*(t_{i-1})$ found in the preceding steps. Hence, equations (5.2) can be solved iteratively using the iterative middle point method [cf. Atkinson (1989)], which then gives approximate values for the unknown boundary $b$ at the points $t_1, \ldots, t_n$.

We recall that the local consistency error for (5.1) for a generic discretization method of the integral is [cf. Linz (1985)]

$$\begin{aligned}(5.3) \quad \delta(h, t_i) &= \int_0^{t_i} \Psi\left(\frac{b(t_i) - b(t_s)}{\sqrt{t_i - s}}\right) f_b(s)\,ds \\ &\quad - h \sum_{j=0}^i \omega_{ij} \Psi\left(\frac{b(t_i) - b(t_j)}{\sqrt{t_i - t_j}}\right) f_b(t_j),\end{aligned}$$

where $\omega_{ij}$ are the integration weights of the adopted discretization schema. In (5.2) we use the Euler method for which we have $\omega_{ij} = 1$ and $\omega_{i0} = 0$, for each $j = 1, \ldots, i$, $i = 1, \ldots, n$. Since $\max_{0 < i < n} \delta(h, t_i) = O(h)$, the method is consistent of order 1.

REMARK 5.1. Note that the system (5.2) is triangular and this makes it especially efficient if one wishes to compute $b^*$ at the next knot, when it is known in the previous ones. Moreover, unlike the PLMC method, here the knowledge of $b(0)$ is not required since we use the "forward" Euler method, also known as right-hand rectangular rule [cf. Atkinson (1989)].

REMARK 5.2. In spite of the Euler method, we could use the extended trapezoidal formula [cf. Abramowitz and Stegun (1964), page 885, formula 25.4.2] with weights $\omega_{i0} = \omega_{ii} = 1/2$ for each $i = 1, \ldots, n$ and $\omega_{ij} = 1$ for each $j = 1, \ldots, i-1$, $i = 1, \ldots, n$, to approximate the integral in (5.1). In this case (5.2) is replaced by

$$
\begin{aligned}
\Psi\left(\frac{b^*(t_i)}{\sqrt{t_i}}\right) = &\frac{1}{2}\Psi\left(\frac{b^*(t_i) - b^*(0)}{\sqrt{t_i}}\right)f_b(0)h \\
&+ \sum_{j=1}^{i-1}\Psi\left(\frac{b^*(t_i) - b^*(t_j)}{\sqrt{t_i - t_j}}\right)f_b(t_j)h \\
&+ \frac{1}{4}f_b(t_i)h \qquad (i = 1, \ldots, n).
\end{aligned}
$$

(5.4)

This approximation is consistent of order 2, but in the general case requests the knowledge of $b(0)$. However, when $f_b(0) = 0$ the first term in the r.h.s. of (5.4) vanishes and we obtain again a triangular system independent from the knowledge of $b(0)$.

REMARK 5.3. The unnecessity of the knowledge of $b(0)$ is a numerical advantage, but it hides some problems connected with the initial point $t_1$. To illustrate this, note that the algorithm uses the approximation (5.2) that for $i = 1$ reads

$$
\Psi\left(\frac{b(t_1)}{\sqrt{t_1}}\right) = \frac{1}{2}f_b(t_1)t_1,
$$
(5.5)

using that $\Psi(0) = 1/2$ (when $b$ is smooth). Recall that $t_1 = t/n$ so that $t_1 \to 0$ as $n \to \infty$. Taking then $b(t) = \sqrt{2t\log(1/t)}$, for example, it is easily verified that $\lim_{n\to\infty} 2\Psi(b(t_1)/\sqrt{t_1})/t_1 = 0$, while the result of Peskir (2002) implies that $f_b(0+) = +\infty$. Thus, the approximation (5.5) in this case fails in a neighborhood of zero. A closer look shows that it is valid if $f_b(0+) = 0$ (under certain mild regularity conditions). It follows that the more $f_b(0+)$ is away from zero, the less accurate (5.5) becomes in a neighborhood of zero.

REMARK 5.4. To improve the boundary estimation for small times, we can consider the flux equation introduced in Ricciardi, Sacerdote and Sato (1984) and its approximation as $t \to 0$ to get

$$(5.6) \qquad f_b(t) \approx \frac{b(t)}{\sqrt{2\pi t^3}} \exp\left[-\frac{b(t)^2}{2t}\right], \qquad t < \epsilon,$$

where $\epsilon$ is small enough to make it acceptable to approximate the r.h.s. of the flux equation with its first term. It is an implicit nonlinear equation in terms of $b(t)$ that can be numerically solved. Since this equation loses reliability as the time $t$ increases, it is necessary to control the validity of this equation for the value of $t$ of interest. Substituting the solution $\hat{b}(t)$ of (5.6) in the complete equation (2.8) in Ricciardi, Sacerdote and Sato (1984), we obtain an estimator $\hat{f}_b(t)$ of the FPT density function. The evaluation of the relative error between this estimator and the known value of the density function allows to guess the time interval where this approximation is reliable.

**6. VIE method error estimation.** Nonlinear integral equations are considered in Linz (1985), but the integral equation (5.1) differs from nonlinear integral equations studied there due to the expression on the l.h.s. Indeed, it contains the unknown function in an implicit way through a further nonlinear function. However, working in analogy with the methods used in Linz (1985), we can prove the convergence of the VIE method of Section 5. To this end, we recall a convergence theorem used by Linz (1985).

THEOREM 6.1. *Let the sequence $\xi_0, \xi_1, \dots$ satisfy*

$$(6.1) \qquad |\xi_n| \leq A \sum_{i=0}^{n-1} |\xi_i| + B_n, \qquad n = r, r+1, \dots,$$

*where $A > 0$, $|B_n| \leq B$ and assume that it exists a constant $\eta > 0$ such that*

$$(6.2) \qquad \sum_{i=0}^{r-1} |\xi_i| \leq \eta.$$

*Then*

$$(6.3) \qquad |\xi_n| \leq (1+A)^{n-r}(B + A\eta), \qquad n = r, r+1, \dots.$$

We now can prove the convergence of the proposed method.

THEOREM 6.2. *The error $\epsilon_n^{\mathrm{VIE}}$ of the VIE method at the discretization knots $t_n$, $n = 1, 2, \dots$, is*

$$(6.4) \qquad |\epsilon_n^{\mathrm{VIE}}| = |b(t_n) - b^*(t_n)| \sim O(h)$$

*if the integral in (5.1) is approximated via the Euler formula.*

*If $f_b(0) = 0$ and the extended trapezoidal formula is used for the integral in (5.1), the error becomes*

$$(6.5) \qquad |\epsilon_n^{\text{VIE}}| = |b(t_n) - b^*(t_n)| \sim O(h^2).$$

PROOF. We split the proof of (6.4) in two parts, the first one concerning the first discretization interval and the second one concerning a generic step $n$.

STEP 1. Choosing $t = t_1$ and subtracting (5.2) from (5.1), we get

$$(6.6) \qquad \Psi\left(\frac{b(t_1)}{\sqrt{t_1}}\right) = \Psi\left(\frac{b^*(t_1)}{\sqrt{t_1}}\right) + \delta(h, t_1).$$

Using the inverse function $\Psi^{-1}$, we obtain

$$(6.7) \qquad \frac{b(t_1)}{\sqrt{t_1}} = \Psi^{-1}\left[\Psi\left(\frac{b^*(t_1)}{\sqrt{t_1}}\right) + \delta(h, t_1)\right].$$

Substituting the Taylor expansion of $\Psi^{-1}$ around $y = \Psi(\frac{b^*(t_1)}{\sqrt{t_1}})$ in this last equation, we get

$$(6.8) \qquad \begin{aligned} \frac{b(t_1)}{\sqrt{t_1}} &= \frac{b^*(t_1)}{\sqrt{t_1}} + \delta(h, t_1)\frac{d\Psi^{-1}(y)}{dy}\bigg|_{y = \Psi(b^*(t_1)/\sqrt{t_1}) + \theta\delta(h, t_1)} \\ &= \frac{b^*(t_1)}{\sqrt{t_1}} + \sqrt{2\pi}\delta(h, t_1)e^{\hat{z}^2/2}, \end{aligned}$$

where $\theta \in (0, 1)$ [cf. Abramowitz and Stegun (1964)] and

$$(6.9) \qquad \hat{z} = \Psi^{-1}\left[\Psi\left(\frac{b^*(t_1)}{\sqrt{t_1}}\right) + \theta\delta(h, t_1)\right] < L$$

with $L < \infty$.

To prove (6.9), note that if $\delta(h, t_1) \geq 0$, one has

$$(6.10) \qquad \Psi\left(\frac{b^*(t_1)}{\sqrt{t_1}}\right) + \theta\delta(h, t_1) \geq \Psi\left(\frac{b^*(t_1)}{\sqrt{t_1}}\right)$$

and, due to the decreasing monotonicity of $\Psi^{-1}$, one gets

$$(6.11) \qquad \hat{z} \leq \frac{b^*(t_1)}{\sqrt{t_1}} < L_1.$$

Elsewhere, when $\delta(h, t_1) \leq 0$ one has

$$(6.12) \quad \Psi\left(\frac{b^*(t_1)}{\sqrt{t_1}}\right) + \theta\delta(h, t_1) \geq \Psi\left(\frac{b^*(t_1)}{\sqrt{t_1}}\right) + \delta(h, t_1) = \Psi\left(\frac{b(t_1)}{\sqrt{t_1}}\right),$$

hence,

$$(6.13) \qquad \hat{z} \le \frac{b(t_1)}{\sqrt{t_1}} < L_2.$$

Choosing $L = \max(L_1, L_2)$, we get (6.9). Hence from (6.8) and (6.9) we get

$$|\epsilon_1^{\mathrm{VIE}}| = |b(t_1) - b^*(t_1)| \le |\sqrt{2\pi t_1}\delta(h, t_1)e^{L^2/2}| < k|\delta(h, t_1)| = O(h),$$

since the Euler method is consistent of order 1.

STEP $n$.  Choosing $t = t_n$ and subtracting (5.2) from (5.1), we get

$$(6.14) \qquad \Psi\left(\frac{b(t_n)}{\sqrt{t_n}}\right) = \Psi\left(\frac{b^*(t_n)}{\sqrt{t_n}}\right) + \gamma_n, \qquad n = 2, 3, \ldots,$$

where

$$
\begin{aligned}
(6.15) \quad & \gamma_n = \gamma_n(h, t_1, \ldots, t_n) \\
& = h\sum_{j=1}^{n} f_b(t_j)\left\{\Psi\left(\frac{b(t_n) - b(t_j)}{\sqrt{t_n - t_j}}\right) - \Psi\left(\frac{b^*(t_n) - b^*(t_j)}{\sqrt{t_n - t_j}}\right)\right\} \\
& \quad + \delta(h, t_n), \qquad n = 2, 3, \ldots.
\end{aligned}
$$

Mimicking the procedure used for case $n = 1$, we apply the inverse function $\Psi^{-1}$ to (6.14) and we expand in Taylor series the resulting r.h.s. around $y = \Psi(\frac{b^*(t_1)}{\sqrt{t_1}})$. Thus, we get

$$
\begin{aligned}
(6.16) \quad \frac{b(t_n)}{\sqrt{t_n}} & = \frac{b^*(t_n)}{\sqrt{t_n}} + \gamma_n\frac{d\Psi^{-1}(y)}{dy}\bigg|_{y = \Psi(b^*(t_n)/\sqrt{t_n}) + \theta\gamma_n} \\
& = \frac{b^*(t_n)}{\sqrt{t_n}} + \sqrt{2\pi}\gamma_n e^{\hat{v}_n^2/2},
\end{aligned}
$$

where $\theta \in (0, 1)$ [cf. Abramowitz and Stegun (1964)] and

$$(6.17) \qquad \hat{v}_n = \Psi^{-1}\left[\Psi\left(\frac{b^*(t_n)}{\sqrt{t_n}}\right) + \theta\gamma_n\right] < M_n$$

with $M_n < \infty$ for each $n = 1, 2, \ldots$.

The last inequality can be proved analogously to (6.11) and (6.13) in Step 1. Hence, from (6.16) we get

$$(6.18) \qquad |\epsilon_n^{\mathrm{VIE}}| = |b(t_n) - b^*(t_n)| \le |\sqrt{2\pi t_n}\gamma_n e^{M_n^2/2}| < k|\gamma_n|.$$

Let us now split the term $\gamma_n$ as the sum of a contribution of the accumulated error due to the previous steps and a contribution due to the consistency of the numerical approximation of the integral in (5.1). To this
goal, we bound the module of $\gamma_n$ using the triangular inequality to get

$$
\begin{aligned}
(6.19) \quad |\gamma_n| &\le h \sum_{j=1}^{n} f_b(t_j) \left| \Psi\left(\frac{b(t_n) - b(t_j)}{\sqrt{t_n - t_j}}\right) - \Psi\left(\frac{b^*(t_n) - b^*(t_j)}{\sqrt{t_n - t_j}}\right) \right| \\
&\quad + |\delta(h, t_n)| \\
&= h \sum_{j=1}^{n-1} f_b(t_j) \left| \frac{1}{\sqrt{2\pi(t_n - t_j)}} \int_{b^*(t_n)-b^*(t_j)}^{b(t_n)-b(t_j)} e^{-u^2/(2(t_n - t_j))} \, du \right| \\
&\quad + |\delta(h, t_n)|,
\end{aligned}
$$

where the last equality uses the definition of $\Psi$ and the fact that the last
term in the sum is zero. Bounding the integrand in (6.19) with 1, we get

$$
(6.20) \quad |\gamma_n| < h \sum_{j=1}^{n-1} \frac{f_b(t_j)}{\sqrt{2\pi(t_n - t_j)}} (|\epsilon_n^{\mathrm{VIE}}| - |\epsilon_j^{\mathrm{VIE}}|) + |\delta(h, t_n)|
$$

and (6.18) becomes

$$
\begin{aligned}
(6.21) \quad |\epsilon_n^{\mathrm{VIE}}| &\le \sqrt{2\pi t_n} M h \sum_{j=1}^{n-1} \frac{f_b(t_j)}{\sqrt{2\pi(t_n - t_j)}} |\epsilon_n^{\mathrm{VIE}}| \\
&\quad + \sqrt{2\pi t_n} M h \sum_{j=1}^{n-1} \frac{f_b(t_j)}{\sqrt{2\pi(t_n - t_j)}} |\epsilon_j^{\mathrm{VIE}}| + \sqrt{2\pi t_n} M |\delta(h, t_n)|,
\end{aligned}
$$

where $M = \max_n \{\exp(M_n)\}$ and, hence,

$$
\begin{aligned}
(6.22) \quad |\epsilon_n^{\mathrm{VIE}}| &\left(1 - M h \sum_{j=1}^{n-1} \frac{f_b(t_j)}{\sqrt{t_n - t_j}}\right) \\
&\le M h \sum_{j=1}^{n-1} \frac{f_b(t_j)}{\sqrt{t_n - t_j}} |\epsilon_j^{\mathrm{VIE}}| + \sqrt{2\pi} M |\delta(h, t_n)|.
\end{aligned}
$$

If

$$
(6.23) \quad h < \frac{1}{M \sum_{j=1}^{n-1} f_b(t_j)/\sqrt{t_n - t_j}} = \frac{1}{k},
$$

we get

$$
(6.24) \quad |\epsilon_n^{\mathrm{VIE}}| \le \frac{MhF}{1 - hk} \sum_{j=1}^{n-1} |\epsilon_j^{\mathrm{VIE}}| + \frac{\sqrt{2\pi} M}{1 - hk} |\delta(h, t_n)|,
$$

where $F = \sup_{j,n} \frac{f_b(t_j)}{\sqrt{2\pi(t_n-t_j)}}$.

Applying the convergence Theorem 6.1 with

$$(6.25) \qquad A_n = \frac{MhF}{1-hk} < A,$$

$$(6.26) \qquad B_n = \frac{2\pi M}{1-hk}|\delta(h,t_n)| < B$$

and noting that by inductive argument one has

$$(6.27) \qquad \sum_{j=1}^{n-1} |\epsilon_j^{\text{VIE}}| \le \eta,$$

we finally obtain

$$(6.28) \quad |\epsilon_n^{\text{VIE}}| \le \left(\max|\delta(h,t_n)|2\pi M + MhF \sum_{j=1}^{n-1} |\epsilon_j^{\text{VIE}}|\right) \frac{1}{1-hk} e^{t_n MF/(1-hk)}.$$

Recalling that the method is consistent of order 1, part 1 of the theorem is proved. Using equation (6.21) in the case of the extended trapezoidal scheme, one easily proves formula (6.5).   $\square$

REMARK 6.1.   Note that the error of the method is dominated by the consistency error.

REMARK 6.2.   To prove Theorem 6.2, we simply use the monotonicity properties of the function $\Psi$. Hence, the method can be easily extended to different diffusion processes, simply substituting $\Psi$ in (5.1) with the survival function of the considered process.

**7. Examples.** In this section we check the stability of the algorithms presented in Sections 3 and 5 by means of some examples where a closed form solution is available. We also show other examples where the solution is numerically evaluated. First we apply the algorithms of Sections 3 and 5 for a Daniels boundary $b$ for two sets of parameters (Section 2). Later we apply them to the case of an oscillating boundary with different parameters. In these later cases the FPT density function has been numerically estimated via the Buonocore, Nobile and Ricciardi algorithm (1987). We consider the mean square deviation as an index of the goodness of the two methods

$$(7.1) \qquad \sigma_n^{(i)} = \frac{1}{n} \sum_{j=0}^{n} (b(t_j) - \hat{b}^{(i)}(t_j))^2, \qquad i = 1,2,$$

where $\hat{b}^{(1)}$ denotes the approximating boundary determined by the PLMC algorithm of Section 3 and $\hat{b}^{(2)}$ denotes the approximating boundary determined by the VIE algorithm of Section 5.

We apply the two algorithms to the following two cases:

1. Daniels boundary with parameters $\alpha = 1, \beta = 0.5, \gamma = 0.5$ and $\alpha = 1, \beta = 1, \gamma = 0.5$. We first apply the PLMC algorithm in the interval $[0, 2]$ with a discretization step $h = 0.2$ and we compute the integrals by a Monte Carlo method using $10^4$ simulations. Under these conditions we estimate the mean square deviation, obtaining $\sigma_n^{(1)} = 9.4 \cdot 10^{-5}$ and $\sigma_n^{(1)} = 3.4 \cdot 10^{-5}$ respectively.

We also apply to these examples the VIE algorithm of Section 5 with discretization step $h = 0.01$ and the resulting mean square deviations are $\sigma_n^{(2)} = 4.3 \cdot 10^{-5}$ and $\sigma_n^{(2)} = 4.6 \cdot 10^{-5}$ respectively.

2. Oscillating boundary

$$(7.2) \qquad b(t) = \alpha + \beta \cos(\gamma t)$$

with parameters $\alpha = 1, \beta = 0.5, \gamma = 2$ and $\alpha = 1, \beta = 1, \gamma = 2$ respectively. Repeating the PLMC algorithm under exactly the same conditions as above, we get the following mean square deviations: $\sigma_n^{(1)} = 6.6 \cdot 10^{-4}$ and $\sigma_n^{(1)} = 6.4 \cdot 10^{-3}$ respectively. The application of the VIE algorithm to the two sets of parameters gives $\sigma_n^{(2)} = 4.9 \cdot 10^{-3}$ and $\sigma_n^{(2)} = 3.4 \cdot 10^{-3}$.

These results confirm the reliability of the methods. The four cases are illustrated in Figures 2, 3, 5 and 6 where the exact boundary shape (left) is compared with the approximating ones obtained by means of the algorithms of Section 3 (center) and Section 5 (right).

REMARK 7.1. With reference to the PLMC method, we underline that the goodness of the approximation does not depend only on the discretization step $h$, but also on the probability mass for the boundary to be crossed in the discretization subinterval of length $h$. It is thus recommendable to avoid excessively small discretization steps. However, as shown in our examples, quite large values of $h$ give small value for the mean square deviation.

REMARK 7.2. A comparison between the exact and the approximating boundary obtained from the PLMC algorithm in Figures 2 and 3 with the first-passage density function in Figure 4 shows the rise of larger oscillations as $t$ increases. This fact can be explained by observing that as $t$ increases the probability of crossing the boundary on each discretized interval becomes smaller and the Monte Carlo method is subject to a larger relative error. A further improvement of the method could be to use an adaptive step built on a constant probability mass of the first-passage density function on each interval.

REMARK 7.3. In Figures 2 and 3 we observe that the approximating boundary obtained from the VIE algorithm of Section 5 has a large error for short times. This result can be explained by the rough approximation in the first steps of the present method and is related to the use of the Euler method. In Remark 5.4 we introduced an improvement of the method for small times. In Figure 7 we compare the approximation obtained via the VIE method (dash dot line) and two different corrections via (5.6) in a special case of Daniels boundary (full line). The stars indicate the values of the boundary solution of (5.6), while the two other curves substitute the first one or two values of $b(t)$ with the solution of (5.6). We observe that the new boundary estimates are more reliable.

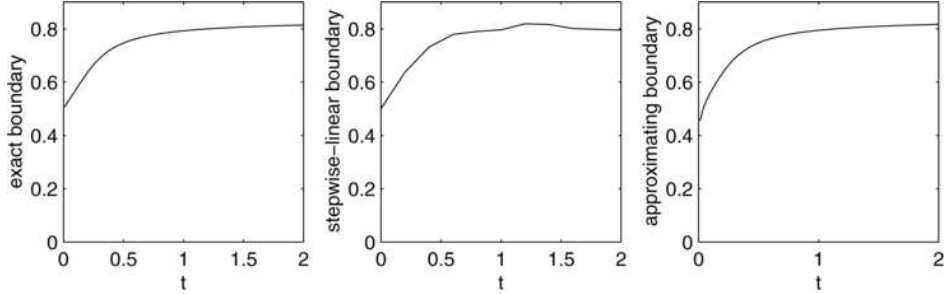

FIG. 2. *Daniels boundary with parameters $\alpha = 1, \beta = 0.5, \gamma = 0.5$ (left) compared with the approximating ones obtained by means of the PLMC method (center) and the VIE method (right).*

REMARK 7.4. The evaluation of the difference between the exact boundary and the approximating one, obtained with the PLMC method, shows that the approximating boundary oscillates around the real one.

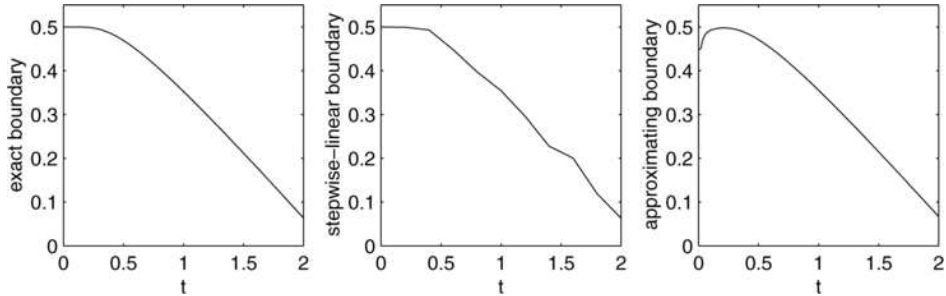

FIG. 3. *Daniels boundary with parameters $\alpha = 1, \beta = 1, \gamma = 0.5$ (left) compared with the approximating ones obtained by means of the PLMC method (center) and the VIE method (right).*

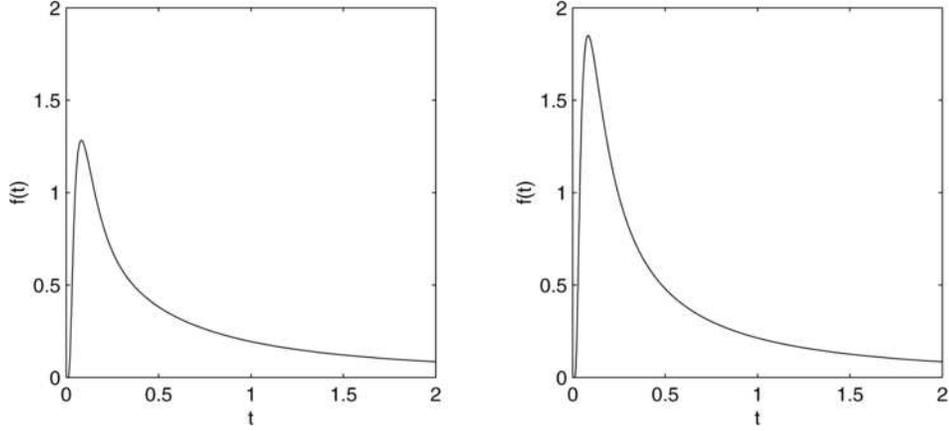

FIG. 4.  *The first-passage density function through a Daniels boundary with parameters*
$\alpha = 1, \beta = 0.5, \gamma = 0.5$ *(left) and with parameters* $\alpha = 1, \beta = 1, \gamma = 0.5$ *(right).*

**8. The exponential case.** Let us apply our algorithms to numerically
solve the original Shiryaev's problem, that is, to approximate the boundary
corresponding to an exponential first-passage density function.

Let us thus consider a first-passage density function $f_b(t) = e^{-t}$ for $t >
0$ corresponding to the exponential distribution with parameter $\lambda = 1$. In
order to apply the PLMC method, we note that this particular choice of
distribution requires to study a first-passage time $\tau_b$ such that $f_b(0+) > 0$.
As indicated in Section 2.2, this implies that the boundary $b$ should be an
upper function for $W$ that vanishes at zero so that $b'(0+) = +\infty$.

Making use of (2.14) and introducing more generally the following nota-
tion:

$$(8.1) \qquad\qquad \kappa = f_b(0+) > 0,$$

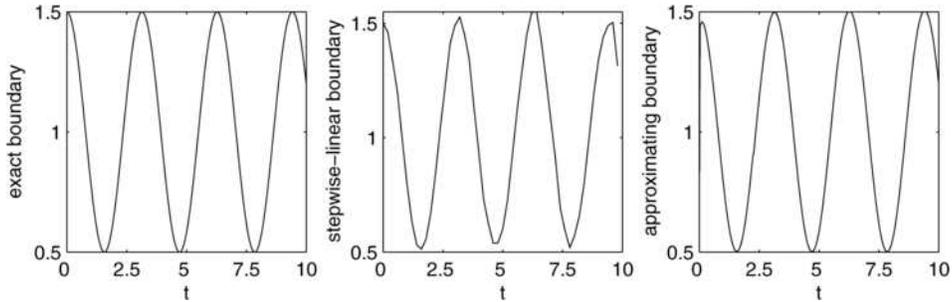

FIG. 5.  *Oscillating boundary with parameters* $\alpha = 1, \beta = 0.5, \gamma = 2$ *(left) compared with
the approximating ones obtained by means of the PLMC method (center) and the VIE
method (right).*

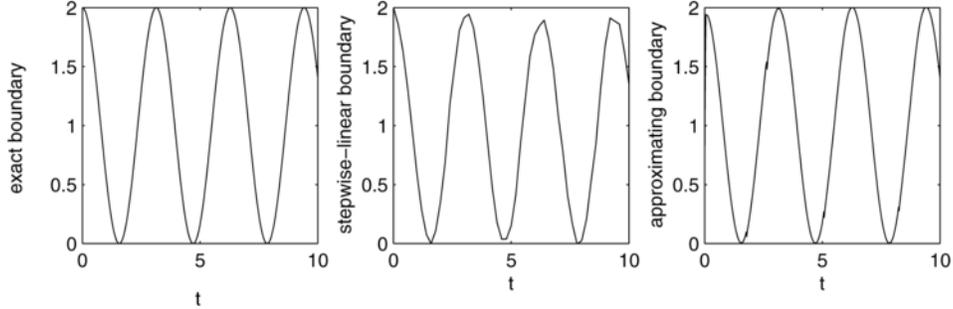

FIG. 6. *Oscillating boundary with parameters $\alpha = 1, \beta = 1, \gamma = 2$ (left) compared with the approximating ones obtained by means of the PLMC method (center) and the VIE method (right).*

we obtain an approximation of the boundary $b$ in the neighborhood of zero by choosing

$$(8.2) \qquad c = -2\log(\sqrt{4\pi}\kappa).$$

On the first interval we approximate the boundary $b$ by (2.13), where $c$ is given by (8.2). On the successive intervals we simply apply the algorithm with

$$(8.3) \qquad c_2 = g(t_1) = \sqrt{2t_1\log(1/t_1) + t_1\log\log(1/t_1) + ct_1}.$$

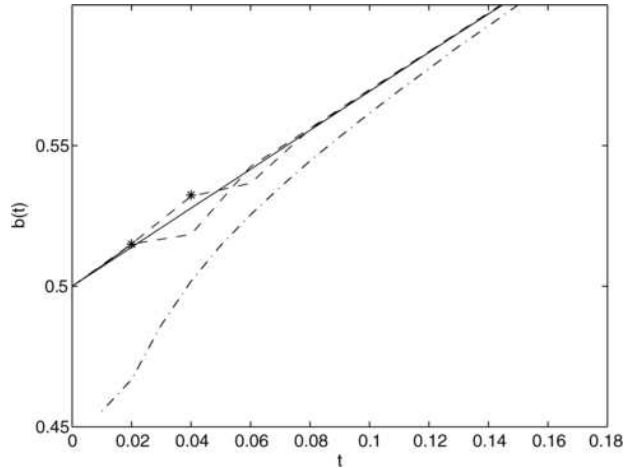

FIG. 7. *Daniels boundary with parameters $\alpha = 1, \beta = 0.5, \gamma = 0.5$ (full line) compared with the approximating ones obtained by means of the VIE method (dash dot line) and by means of the VIE method modified for small times (dashed lines) substituting the first one or two values of $b(t)$ with the solution of (5.6). The stars indicate the values of the boundary obtained as a solution of (5.6).*

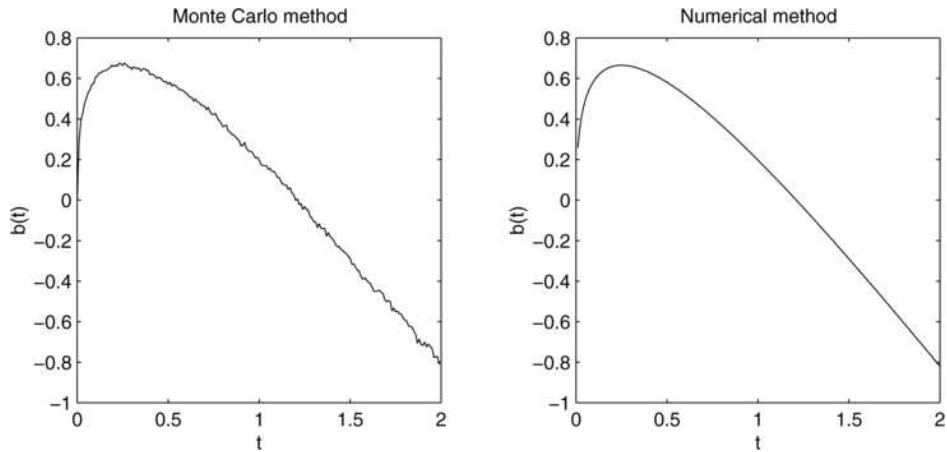

FIG. 8. *Numerical evaluations of the approximating boundaries for an exponential first–passage density function with parameter $\lambda = 1$ obtained by means of the PLMC method (left) and the VIE method (right) with a discretization step $h = 0.01$.*

In Figure 8 we plot the boundaries corresponding to an exponential first-passage density function with $\lambda = 1$ obtained by the two different algorithms described in the previous sections; the plot on the r.h.s. is obtained applying the PLMC method with the step $h = 0.01$ and with a number of simulations for the Monte Carlo method equal to $10^4$, while the plot on the l.h.s. is obtained applying the VIE method with discretization step $h = 0.01$.

A heuristic confirmation of the stability of the algorithms is given by the nearness of the two curves. Furthermore, an intuitive reasoning about the shape of the boundary confirms our result. To obtain an exponential first-passage density function, a large part of the sample paths should cross the boundary very soon but, as the time goes on, other sample paths, that have not yet reached the boundary, should be able to cross it too. This makes intuitive the fact that the boundary must decrease in order to be reachable by the sample paths whose initial trend was negative.

**Acknowledgments.** We are grateful to Professor G. Peskir for bringing the problem to our attention and for his interesting and useful comments in the first development of this paper. We also acknowledge Dr. W. Dambrosio for useful suggestions and two anonymous referees for their help to improve the paper.

## REFERENCES

ABRAMOWITZ, M.and STEGUN, I. A. (1964). *Handbook of Mathematical Functions With Formulas, Graphs, and Mathematical Tables.* Dover, New York.

ANULOVA, S. V. (1980). On Markov stopping times with a given distribution for a Wiener process. *Theory Probab. Appl.* **5** 362–366.

ATKINSON, K. E. (1989). *An Introduction to Numerical Analysis*, 2nd ed. Wiley, New York. MR1007135

BUONOCORE, A., NOBILE, A. G. and RICCIARDI, L. M. (1987). A new integral equation for the evaluation of first-passage-time probability densities. *Adv. in Appl. Probab.* **19** 784–800. MR914593

CAPOCELLI, R. M. and RICCIARDI, L. M. (1972). On the inverse of the first passage time probability problem. *J. Appl. Probability* **9** 270–287. MR0346926

DANIELS, H. E. (1969). The minimum of a stationary Markov process superimposed on a *U*-shaped trend. *J. Appl. Probability* **6** 399–408. MR0264782

DOOB, J. L. (1949). Heuristic approach to the Kolmogorov–Smirnov theorems. *Ann. Math. Statist.* **20** 393–403. MR0030732

DUDLEY, R. M. and GUTMANN, S. (1977). Stopping times with given laws. In *Séminaire de Probabilités, XI (Univ. Strasbourg, Strasbourg, 1975/1976). Lecture Notes in Math.* **581** 51–58. Springer, Berlin. MR0651553

DURBIN, J. (1971). Boundary-crossing probabilities for the Brownian motion and Poisson processes and techniques for computing the power of the Kolmogorov–Smirnov test. *J. Appl. Probab.* **8** 431–453. MR0292161

LINZ, P. (1985). *Analytical and Numerical Methods for Volterra Equations. SIAM Studies in Applied Mathematics* **7**. SIAM, Philadelphia, PA. MR796318

MALMQUIST, S. (1954). On certain confidence contours for distribution functions. *Ann. Math. Statistics* **25** 523–533. MR0065095

PESKIR, G. (2002). Limit at zero of the Brownian first-passage density. *Probab. Theory Related Fields* **124** 100–111. MR1929813

RICCIARDI, L. M., SACERDOTE, L. and SATO, S. (1984). On an integral equation for first-passage-time probability densities. *J. Appl. Probab.* **21** 302–314. MR741132

RICCIARDI, L. M. and SATO, S. (1990). Diffusion processes and first-passage-time problems. In *Lectures in Applied Mathematics and Informatics* 206–285. Manchester Univ. Press, Manchester. MR1075228

SACERDOTE, L. and SMITH, C. E. (2004). Almost sure comparisons for first passage times of diffusion processes through boundaries. *Methodol. Comput. Appl. Probab.* **6** 323–341. MR2073545

SACERDOTE, L., VILLA, A. E. P. and ZUCCA, C. (2006). On the classification of experimental data modeled via a stochastic leaky integrate and fire model through boundary values. *Bull. Math. Biol.* **68** 1257–1274. MR2249750

SACERDOTE, L. and ZUCCA, C. (2003a). Threshold shape corresponding to a gamma firing distribution in an Ornstein–Uhlenbeck neuronal model. *Sci. Math. Jpn.* **58** 295–305. MR2021074

SACERDOTE, L. and ZUCCA, C. (2003b). On the relationship between interspikes interval distribution and boundary shape in the Ornstein–Uhlenbeck neuronal model. In *Mathematical Modelling & Computing in Biology and Medicine. Milan Research Centre for Industrial and Applied Mathematics. The MIRIAM Project* **1** 161–167. Esculapio, Bologna. MR2093270

STRASSEN, V. (1967). Almost sure behavior of sums of independent random variables and martingales. In *Proceedings Fifth Berkeley Symposium Math. Statist. Probab. Vol. II: Contributions to Probability Theory, Part 1* 315–343. Univ. California Press, Berkeley. MR0214118

WANG, L. and PÖTZELBERGER, K. (1997). Boundary crossing probability for Brownian motion and general boundaries. *J. Appl. Probab.* **34** 54–65. MR1429054

Zucca, C. (2002). Analytical, numerical and Monte Carlo techniques for the study of the first passage times. Ph.D. thesis, Univ. Milano.

Department of Mathematics
University of Torino
Via Carlo Alberto 10
10123 Torino
Italy
E-mail: cristina.zucca@unito.it
           laura.sacerdote@unito.it
URL: www2.dm.unito.it/paginepersonali/zucca
           www.dm.unito.it/personalpages/sacerdote/index.htm